\documentclass[11pt]{amsart}
\usepackage{amsfonts}
\usepackage{amssymb,latexsym}
\usepackage{enumerate}
\usepackage{mathrsfs}
\usepackage{amsmath}
\usepackage{amsthm}
\newtheorem{Thm}{Theorem}[section]

\newtheorem{Rem}[Thm]{Remark}

\numberwithin{equation}{section}

\begin{document}
\setlength{\baselineskip}{1.2\baselineskip}
\title[Gradient Estimates]{Gradient Estimates of Mean Curvature Equation with Neumann Boundary Condition in $ M^n\times \mathbb{R}$}
\author{Jinju Xu}
\address{Department of Mathematics\\
         Shanghai University\\
         Shanghai 200444 CHINA}
         \email{jjxujane@shu.edu.cn}
\author{Dekai Zhang}
\address{School of Mathematical Sciences\\
         University of Science and Technology of China\\
         Hefei Anhui 230026 CHINA}
         \email{dekzhang@mail.ustc.edu.cn}


\maketitle
\begin{abstract}
 In this note, we study the prescribed mean curvature equation with
 Neumann boundary conditions on Riemannian product manifold $ M^n\times \mathbb{R}$. The main goal is to establish  the  boundary gradient estimates for solutions by the maximum principle.  As a consequence, we obtain an existence result.
\end{abstract}

\section{Introduction }
In this paper,  we consider the following Neumann problem
\begin{align}
\texttt{div}(\frac{Du}{\sqrt{1+|Du|^2}})=& f(x,u)\quad\text{in}\quad\Omega,\label{1.1}\\
<Du,\gamma>=&\psi(x,u)\quad\text{on}\quad\partial\Omega\label{nu},
\end{align}
where  $\Omega$ is a bounded domain in $n$-dimensional  simply connected and complete manifold $ M\subset\mathbb{R}^{n+1}$ with Riemannian metric $\sigma$, $f$ and $\psi$ are given functions on $M\times \mathbb{R}$ and $\partial\Omega\times \mathbb{R}$ respectively,  $\gamma$ is the unit
inner normal to $\partial\Omega\times \mathbb{R}$, and $<,>$ denotes the Riemannian metric in $M$.

Let $S=\{\left(x, u(x)\right):x\in \Omega\}$ be  $n$ dimensional graphs of mean curvature $H$ in an $n+1$  dimensional Riemannian manifold of the form $ M^n\times \mathbb{R}$.
  If $ds^2=\sigma_{ij}dx_idx_j$ is a local Riemannian metric on $ M $,
 then $M\times \mathbb{R}$ is given the product metric $ds^2 +dt^2$ where $t$ is a coordinate for $\mathbb{R}$. Then the height function $u(x)\in C^2(\Omega)$ satisfies the following equation
 \begin{align}
 \texttt{div}(\frac{Du}{\sqrt{1+|Du|^2}})=nH(x),\label{1.2}
 \end{align}
where the divergence and gradient $Du$ are taken with respect to the metric on $M$ i.e.
\begin{align*}
Du: = \sum\limits_{i = 1}^n {u^i \frac{\partial }{ {\partial x_i }}}, \quad u^i  = \sum\limits_{j = 1}^n {\sigma ^{ij} \frac{\partial u}{\partial x_j }}.
\end{align*}
The equation \eqref{1.1} is equivalent to the following
\begin{align}\label{equ}
 \sum\limits_{i,j = 1}^ng^{ij}D_iD_ju=&f(x,u)v\quad\text{in}\quad\Omega,\quad
  v=\sqrt{1+|Du|^2},
\end{align}
where $D$ denotes covariant differentiation on $M$ and
\begin{align*}
g^{ij}=\sigma^{ij}-\frac{u^i u^j}{1+|Du|^2}.
\end{align*}

For the following prescribed mean curvature equation with prescribed contact angle boundary value problem in $\mathbb{R}^{n+1}(n\ge2)$,
\begin{equation}
\left\{\begin{array}{rcl}\texttt{div}(\frac{Du}{\sqrt{1+|Du|^2}})=& f(x,u)&\quad\text{in}\quad\Omega,\\
<Du,\gamma>=&\phi (x)v&\quad\text{on}\quad\partial\Omega,
\end{array} \right.
\end{equation}
Ural'tseva \cite{Ur73} first
got the boundary gradient estimates and the corresponding positive gravity case existence theorem. At the same time, Simon-Spruck \cite{SS76} and Gerhardt \cite{Ger76}
also obtained existence theorem on the positive gravity case. They obtained these estimates also via test function technique.
Spruck \cite{Sp75} used the maximum principle to obtain boundary gradient estimate in two dimension for the positive gravity capillary problems.
 Korevaar \cite{Kor88} generalized his normal variation technique and got the gradient estimates for the positive gravity case in high dimension case.
 Lieberman \cite{Lieb88}  got  the gradient estimates using a closely related maximum principle argument to Korevaar's  method on more general quasilinear elliptic equations with capillary boundary value problem in zero gravity case.
Recently in \cite{Xu14}, the first author used the
  maximum principle  to give a new proof of gradient estimates for  the mean curvature equation with oblique problem including capillary boundary value problem in zero gravity case.


In Riemannian product manifolds $ M^n\times \mathbb{R}$,
Spruck \cite{Sp07} proved the interior gradient estimates and existence theorems of Dirichlet problem for  constant mean curvature graphs. Many researchers have also studied the capillary problem, see \cite{LW13} and references therein.
Recently, Ma-Xu \cite{MX14} used the maximum principle to get the gradient estimate
for the solutions of the  mean curvature equation with Neumann boundary value and  obtain an existence result. Naturally, we want to know whether  we  can generalize the Euclidean result  to Riemannian product manifolds.

In this paper, we use Ma-Xu's technique in \cite{MX14} to give the boundary gradient estimates for general mean curvature  equation of graphs with Neumann boundary condition in Riemannian product manifolds $ M^n\times \mathbb{R} $. Here our proof is more simple than that in \cite{MX14}.

Now we state our main result.
\begin{Thm}\label{Thm1.1} Let $\Omega \subset M^n $ be a bounded domain, $n\geq 2$, $\partial \Omega \in C^{3}$, $\gamma$ is the inward unit normal vector to $\partial\Omega $.
Suppose $u\in C^{2}(\overline\Omega)\bigcap C^{3}(\Omega)$ is a bounded solution of $\eqref{1.1}, \eqref{nu}$ with $|u|\le M_0$ where $M_0$ is a positive constant.
Assume $f(x,z) \in C^{1}(\overline\Omega\times [-M_0, M_0])$ and $\psi(x,z) \in C^{3}(\overline\Omega\times [-M_0, M_0])$,
and there exist positive constants $ L_1, L_2$ such that
\begin{align}
f_z(x,z)\geq &0 \quad \text{in}\quad \overline\Omega\times [-M_0, M_0],\label{1.5}\\
|f(x,z)|+|D_{x}f(x,z)|\leq & L_1 \quad \text{in}\quad \overline\Omega\times [-M_0, M_0],\label{1.6}\\
|\psi(x,z)|_{C^3(\overline\Omega\times[-M_0, M_0])}\leq&  L_2.\label{1.7}
\end{align}
Then there exists a small positive constant
$\mu_0$ such that we have the following estimate
$$\sup_{\overline\Omega_{\mu_0}}|Du|\leq \max\{M_1, M_2\},$$
where $M_1$ is a positive constant depending only on $n, \mu_0, M_0, L_1, \sup_{\overline\Omega}|Ric|$, which is from the interior gradient estimates;
$M_2$ is  a positive constant depending only on $n, \Omega, \mu_0, M_0, L_1, L_2, \sup_{\overline\Omega}|Ric|$ and $d(x) =\texttt{dist}(x, \partial\Omega), \Omega_{\mu_0} = \{x \in \Omega: d(x)<\mu_0\}.$
\end{Thm}
As in \cite{Sp07}, there is an interior gradient estimates for the mean curvature equation in $ M^n\times \mathbb{R} $. One can also use the method of Trudinger in \cite{Tru90} or  Xu-jia Wang's method in \cite{Wang98} to give the
interior gradient estimates for more general  mean curvature equations in $ M^n\times \mathbb{R} $. We state it in the following.

\begin{Rem}\label{Rem1.1}
If $u\in C^{3}(\Omega)$ is a bounded solution for the equation \eqref{1.1}  with  $|u|\le M_0$ , and if $f \in C^{1}(\overline\Omega \times [-M_0, M_0])$
satisfies the conditions \eqref{1.5}-\eqref{1.6}, then for any subdomain
$\Omega'\subset\subset\Omega$, we have
$$\sup_{\Omega'}|Du|\leq M_1,$$
where $M_1$ is a positive constant depending only on $n,  M_0, \texttt{dist} (\Omega', \partial\Omega), L_1, \sup_{\overline\Omega}|Ric|$.
\end{Rem}

From the standard bounded estimates for the prescribed mean curvature equation in Concus-Finn \cite{CF74} ( see also Spruck \cite{Sp75}, Ma-Xu\cite{MX14}), we can also get the following existence theorem.

\begin{Thm}\label{Thm1.2}
Let $\Omega \subset M^n $ be a bounded domain, $n\geq 2$, $\partial \Omega \in C^{3}$, $\gamma$ is the inward unit normal vector to $\partial\Omega $.
If  $\psi \in C^{3}(\overline\Omega)$  is a given function, then the following boundary value problem
\begin{equation}
\left\{\begin{array}{rl}\texttt{div}(\frac{Du}{\sqrt{1+|Du|^2}})=u&\quad\text{in}\quad\Omega,\\
<Du,\gamma>=\psi (x)&\quad\text{on}\quad\partial\Omega,
\end{array} \right.
\end{equation}
 exists a unique solution $ u \in C^2(\overline\Omega)$.
\end{Thm}

The rest of the paper is organized as follows. In section 2, we first give  some preliminaries. Then
 we shall prove the main Theorem~\ref{Thm1.1} in section 3. As a corollary, we obtain the existence result Theorem~\ref{Thm1.2}.

\section{Preliminaries}
In this section, we give some notations which are mainly from \cite{Sp07} and \cite{MX14}.
Let $x_1,\ldots,x_n$ be a system of local coordinates for $M$ with corresponding metric $\sigma_{ij}$. Then the coordinate vector fields for $S$ and the upward unit normal to $S$
is given by
\begin{align}\label{2.1}
X_i=&\frac{\partial}{\partial x_i}+u_i\frac{\partial}{\partial t},
\end{align}
and
\begin{align}\label{2.2}
N=&\frac{1}{v}(-\sum_{1\leq i\leq n}u^i\frac{\partial}{\partial x_i}+\frac{\partial}{\partial t}),\quad
u^i=\sum_{1\leq j\leq n}\sigma ^{ij}u_j.
\end{align}
The induced metric on $S$ is then
\begin{align}\label{2.3}
g_{ij}=&<X_i,X_j>=\sigma_{ij}+u_iu_j,
\end{align}
with inverse
\begin{align}\label{2.4}
g^{ij}=&\sigma^{ij}-\frac{u^i u^j}{v^2}.
\end{align}
The second fundamental form $b_{ij}$ of $S$ is given by ($\overline{D}$ is covariant differentiation on $M^n\times \mathbf{R}$)
\begin{align*}
b_{ij}=&<\overline{D}_{X_i}X_j,N>\\
=&<D_{\frac{\partial}{\partial x_i}}\frac{\partial}{\partial x_j}+u_{ij}\frac{\partial}{\partial t},N>\\
=&<\sum_{1\leq k\leq n}\Gamma^k_{ij}\frac{\partial}{\partial x_k}+u_{ij}\frac{\partial}{\partial t},N>\\=&\frac{1}{v}(-\Gamma^k_{ij}\sum_{1\leq l\leq n}u^l\sigma_{kl}+u_{ij}).
\end{align*}
Hence
\begin{align}\label{2.5}
b_{ij}=&\frac{1}{v}D_iD_ju
\end{align}
and the equation of prescribed mean curvature $H(x)$ is then
\begin{align}\label{2.6}
nH(x)=&\frac{1}{v}\sum_{1\leq i,j\leq n}g^{ij}D_iD_ju,
\end{align}
where we use the notation $v=\sqrt{1+|Du|^2}$ for convenience.

We denote by $\Omega$ a bounded  domain in $M^n$, $n\geq 2$,  $\partial \Omega\in C^{3}$,  and  set
\begin{align*}
 d(x)=\texttt{dist}(x,\partial \Omega),
 \end{align*}
 and
\begin{align*}
 \Omega_\mu=&\{{x\in\Omega:d(x)<\mu}\}.
 \end{align*}
 It is well known that there exists a positive constant $\mu_{1}>0$ such that $d(x)\in C^3(\overline \Omega_{\mu_{1}}) ,|Dd|= 1$.

As in \cite{MX14}, we define
 \begin{align}\label{2.9}
c^{ij}=&\sigma ^{ij}-\gamma^i\gamma^j  \quad \text{in} \quad \Omega_{\mu_{1}},
\end{align}
 and for a vector $\zeta \in R^n$, we write $\zeta'$ for the vector with $i-$th component $ \sum_{1\leq j\leq n}c^{ij}\zeta_j$. So
 \begin{align}\label{2.10}
|D'u|^2=& \sum_{1\leq i,j\leq n}c^{ij}u_iu_j.
\end{align}

\section{Proof of Theorem~\ref{Thm1.1} }
Now we begin to prove Theorem~\ref{Thm1.1}, using the technique developed by Ma-Xu \cite{MX14}. Here our calculation is more simple than that in \cite{MX14}. We shall choose an auxiliary function which contains  $|Du|^2$ and other lower order terms. Then we use the maximum principle for this auxiliary function in $\overline\Omega_{\mu_0}, 0<\mu_0<\mu_1$. At last,  we get our estimates.

{\em Proof of  Theorem~\ref{Thm1.1}.}

Setting $w=u-\psi(x,u)d$,  we choose the following auxiliary function
$$\Phi(x)=\log|Dw|^2e^{1+M_0+u}e^{\alpha_0 d}, \quad  x \in \overline \Omega_{\mu_{0}},$$ where $\alpha_0=|\psi|_{C^0(\overline\Omega\times[-M_0, M_0])}+C_0+1 $, $C_0$ is  a positive constant depending only on $n,\Omega$.

In order to simplify the computation, we let
\begin{align}\label{3.0}
 \varphi(x)= \log \Phi(x) =\log\log|Dw|^2+h(u)+g(d),
 \end{align}
where we take
\begin{align}\label{3.1}
 h(u)=1+M_0+u,\quad g(d)=\alpha_0 d.
 \end{align}

 We assume that
$\varphi(x)$ attains its maximum at $x_0 \in \overline \Omega_{\mu_{0}}$, where $0<\mu_0<\mu_1$ is a sufficiently small number which we shall decide it  later.

Now we divide three cases to complete the proof of  Theorem~\ref{Thm1.1}.

Case I. If $\varphi(x)$ attains its maximum at $x_0 \in \partial\Omega$, then we shall  get the bound of $|Du|(x_0)$.

Case II. If $\varphi(x)$ attains its maximum at $x_0 \in\partial\Omega_{\mu_0}\bigcap\Omega$, then we shall
 get the estimates via the  interior gradient bound in Remark~\ref{Rem1.1}.

Case III.  If $\varphi(x)$ attains its maximum at $x_0 \in \Omega_{\mu_0}$,  then we can
use the maximum principle to get the bound of $|Du|(x_0)$.

Now  all computations will be done  at the point $x_0$.

{\bf Case I.} If  $\varphi(x)$ attains its maximum at $x_0\in \partial \Omega$, we shall get the bound of $|Du|(x_0)$.  We choose the normal coordinate at $x_0$, such that
\begin{align*}
\sigma _{ij}(x_0)= \delta _{ij},\quad
 w_1(x_0) =|Dw|(x_0).
\end{align*}

We differentiate $\varphi$ along the normal direction.
\begin{align}\label{3.2}
\frac{\partial\varphi}{\partial\gamma}=&\frac{\sum_{1\leq i\leq n}D_i(|Dw|^2)\gamma^i}{|Dw|^2\log|Dw|^2}+h'u_{\gamma}+g'.
\end{align}
Since
\begin{align}
D_iw=&D_iu-\psi_u D_iu d -D_{x_i}\psi d-\psi\gamma^i,\label{3.3}\\
|Dw|^2=&|D'w|^2+w^2_\gamma,\label{3.4}
\end{align}
we have
\begin{align}
w_\gamma=&u_\gamma-\psi_uu_\gamma d-D_{x_i}\psi\gamma^i d-\psi=0\quad\text{on}\quad\partial\Omega,\label{3.5}\\
D_i(|Dw|^2)=&D_i(|D'w|^2)\quad\text{on}\quad\partial\Omega.\label{3.6}
\end{align}
Applying 
\eqref{2.9} and \eqref{3.6}, it follows that
\begin{align}
\sum_{1\leq i\leq n}D_i(|Dw|^2)\gamma^i=&\sum_{1\leq i\leq n}D_i(|D'w|^2)\gamma^i\notag\\
=&2\sum_{1\leq i,k,l\leq n}c^{kl}D_iD_kuD_lu\gamma^i-2\sum_{1\leq k,l\leq n}c^{kl}D_luD_{k}\psi,\label{3.7}
\end{align}
where
\begin{align*}
D_k\psi= D_{x_k}\psi + \psi_u D_ku.
\end{align*}
Differentiating \eqref{nu} with respect to tangential direction,   we have
\begin{align}\label{3.8}
\sum_{1\leq k\leq n}c^{kl}D_k(u_{\gamma})=&\sum_{1\leq k\leq n}c^{kl}D_k\psi.
\end{align}
It follows that
\begin{align}\label{3.9}
\sum_{1\leq i,k\leq n}c^{kl}D_kD_iu\gamma^i=&-\sum_{1\leq i,k\leq n}c^{kl}D_iuD_k(\gamma^i)+\sum_{1\leq k\leq n}c^{kl}D_k\psi.
\end{align}
Inserting \eqref{3.9} into \eqref{3.7} and combining \eqref{nu}, \eqref{3.2}, we have
\begin{align}\label{3.10}
|Dw|^2\log|Dw|^2\frac{\partial\varphi}{\partial\gamma}(x_0)
=&(g'(0)+h'\psi)|Dw|^2\log|Dw|^2-2\sum_{1\leq i,k,l\leq n}c^{kl}D_iuD_luD_k(\gamma^i).
\end{align}
From \eqref{3.3}, we obtain
\begin{align}\label{3.11}
|Dw|^2=&|Du|^2-\psi^2\quad\text{on}\quad\partial\Omega.
\end{align}
Assume $|Du|(x_0)\ge \sqrt{100+2|\psi|^2_{C^0(\overline\Omega\times[-M_0, M_0])}}$, otherwise we get the estimates. At $x_0$,   we have
\begin{align}
\frac{1}{2}|Du|^2\leq&|Dw|^2\leq |Du|^2\quad\text{and}\quad|Dw|^2\ge 50.\label{3.12}
\end{align}
Inserting  \eqref{3.12}  into \eqref{3.10}, and by the choice of $\alpha_0$, we have
\begin{align}\label{3.14}
\frac{\partial\varphi}{\partial\gamma}(x_0)
\geq&\alpha_0-|\psi|_{C^0(\overline\Omega\times[-M_0, M_0])}-\frac{2\sum_{1\leq i,k,l\leq n}|c^{kl}D_iuD_luD_k(\gamma^i)|}{|Dw|^2\log|Dw|^2}\notag\\
\geq&\alpha_0-|\psi|_{C^0(\overline\Omega\times[-M_0, M_0])}-C_0\notag\\
>&0.
\end{align}
On the other hand,  it is obvious to get
$$\frac{\partial\varphi}{\partial\gamma}(x_0)\leq 0,$$
which is a contradiction to \eqref{3.14}.

Then we have
\begin{align}\label{3.15}
|Du|(x_0)\leq\sqrt{100+2|\psi|^2_{C^0(\overline\Omega\times[-M_0, M_0])}}.
\end{align}

{\bf Case II.} $ x_0\in \partial\Omega_{\mu_{0}}\bigcap\Omega$.
This is due to interior gradient estimates. From Remark~\ref{Rem1.1}, we have
 \begin{align}\label{3.16}
\sup_{\partial\Omega_{\mu_0}\bigcap\Omega}|Du|\leq \tilde{M}_1,
\end{align}
where $\tilde{M}_1$ is a positive constant depending only on $n, M_0, \mu_0, L_1$.\par

{\bf Case III.} $x_0\in\Omega_{\mu_{0}}$. \par

In this case, $x_0$ is a critical point of $\varphi$. We choose the normal coordinate at $x_0$, such that
\begin{align*}
\sigma _{ij}(x_0)= \delta _{ij},\quad
 w_1(x_0) =|Dw|(x_0).
\end{align*}
 And  the matrix $(D_iD_jw(x_0))(2\leq i,j\leq n)$ is diagonal. Let $$\mu_2 \le \frac{1}{10
0L_2}$$ such that
 \begin{align}\label{3.16a}
 |\psi_u| \mu_2 \le \frac{1}{100
}, \quad \text{then}\quad \frac{99}{100}\le 1-\psi_u \mu_2 \le
 \frac{101}{100}.
\end{align}
 We can choose $$ \mu_0=\frac{1}{2} \min\{ \mu_1, \mu_2,  1 \}.$$
In order to simplify the calculations, we let
\begin{align*}
w=&u-G,\quad G=\psi(x,u)d.
\end{align*}
Then we have
\begin{align}
w_k=(1-G_u)u_k-G_{x_k}.\label{3wk}
\end{align}
Since at $x_0$,
\begin{align}
|Du|^2=&u_1^2+\sum_{2\leq i\leq n}u_i^2,\label{3Aa}\\
(1-G_u)u_i=&G_{x_i}=\psi_{x_i}d+\psi d_i,\quad i=2,\ldots,n,\label{3wi2}\\
w_1=&(1-G_u)u_1-G_{x_1}=(1-G_u)u_1 -\psi_{x_1}d-\psi d_1.\label{3wi2a}
\end{align}
So from the above relations,  at $x_0$, we can assume
\begin{align}
u_1(x_0)\ge 200(1+ |\psi|_{C^1(\overline\Omega\times[-M_0, M_0])}),\label{3u1}
\end{align}
then
\begin{align}\label{3.18}
\frac{10}{11}w_1^2\leq|Du|^2\leq \frac{11}{9}w_1^2, \quad
\frac{19}{20}u_1\leq w_1\leq \frac{21}{20}u_1,
\end{align}
and by
 the choice of $\mu_0$ and  \eqref{3.16a}, we have
\begin{align}
\frac{99}{100}
\le 1-G_u \le\frac {101}{100}.
\label{wi2aa}
\end{align}

From the above choice, we shall prove Theorem~\ref{Thm1.1} with two steps. As we mentioned before, all the calculations will be done at the fixed point $x_0$. In the following, we denote by
$D_iu=u_i,\, D_jD_iu=u_{,ij},\,  D_kD_jD_iu=u_{,ijk},\ldots$

{\bf Step 1:} We first get the formula \eqref{3aijvarphiija} and the estimate \eqref{3I2c}.\par

Taking the first covariant derivatives of $\varphi$,
\begin{align}\label{3varphii}
\varphi_i=&\frac{(|Dw|^2)_i}{|Dw|^2\log|Dw|^2}+h'u_i+g'd_i.
\end{align}
From $\varphi_i(x_0)=0$,  we have
\begin{align}\label{3varphii=0}
(|Dw|^2)_i=-|Dw|^2\log|Dw|^2(h'u_i+g'd_i).
\end{align}
Taking the covariant derivatives again for  $\varphi_i$, we have
\begin{align}
\varphi_{,ij}
=&\frac{(|Dw|^2)_{,ij}}{|Dw|^2\log|Dw|^2}-(1+\log|Dw|^2)\frac{(|Dw|^2)_i(|Dw|^2)_j}{(|Dw|^2\log|Dw|^2)^2}\notag\\&
+h'u_{,ij}
+h''u_iu_j+g''d_id_j+g'd_{,ij}.\label{3varphiija}
\end{align}
Using \eqref{3varphii=0}, it follows that
\begin{align}
\varphi_{,ij}
=&\frac{(|Dw|^2)_{,ij}}{|Dw|^2\log|Dw|^2}+h'u_{,ij}+\big[h''-(1+\log|Dw|^2)h'^2\big]u_iu_j+g'd_{,ij}\notag\\&
+\big[g''-(1+\log|Dw|^2)g'^2\big]d_id_j
-(1+\log|Dw|^2)h'g'(d_iu_j+d_ju_i).\label{3varphiijb}
\end{align}
Then we get
\begin{align}\label{3aijvarphiija}
0\geq \sum_{1\leq i,j\leq n}g^{ij}\varphi_{,ij} (x_0)
=: & I_1+I_2,
\end{align}
where
\begin{align}\label{3I1a}
I_1=\frac{1}{|Dw|^2\log|Dw|^2}\sum_{1\leq i,j\leq n}g^{ij}(|Dw|^2)_{,ij},
\end{align}
and
\begin{align}
I_2=&\sum_{1\leq i,j\leq n}g^{ij}\big\{h' u_{,ij}+\big[h''-(1+\log|Dw|^2)h'^2\big]u_iu_j\notag\\
&+\big[g''-(1+\log|Dw|^2)g'^2\big]d_id_j-2(1+\log|Dw|^2)h'g'd_iu_j+g'd_{,ij}\big\}.\label{3I2a}
\end{align}

From the choice of the coordinate and the equations \eqref{equ}, \eqref{3.1}, we have the estimate for $I_2$:
\begin{align}
I_2
\ge&f v-2\log w_1\sum_{1\leq i,j\leq n}g^{ij}u_iu_j-2\alpha_0^2\log w_1\sum_{1\leq i,j\leq n}g^{ij}d_id_j-C_1,\label{3I2c}
\end{align}
here  $C_1$   depending  only on $n, \Omega, M_0, \mu_0, L_2$.\par

{\bf Step 2:} We  treat  $I_1$  and finish the proof of  Theorem~\ref{Thm1.1}.\par

Taking the first covariant derivatives of $|Dw|^2$, we have
\begin{align}\label{3|Dw|^2i}
(|Dw|^2)_i=&2w_1w_{,1i}.
\end{align}
Taking the covariant derivatives of $|Dw|^2$ once more, we have
\begin{align}\label{3|Dw|^2ij}
(|Dw|^2)_{,ij}=&2w_1w_{,1ij}+2\sum_{1\leq k\leq n}w_{,ki}w_{,kj}.
\end{align}
By \eqref{3I1a} and \eqref{3|Dw|^2ij}, we can rewrite $I_1$ as
\begin{align}
I_1
=&\frac{1}{w_1^2\log w_1}\big[\sum_{1\leq i,j\leq n}g^{ij}w_1w_{,1ij}+\sum_{1\leq i,j,k\leq n}g^{ij}w_{,ki}w_{,kj}\big]\notag\\
=:&\frac{1}{w_1^2\log w_1}\big[I_{11}+I_{12}\big].\label{3I1b}
\end{align}
In the following, we shall deal with $I_{11}$ and $I_{12}$ respectively. \par

For the term $I_{11}=w_1\sum_{1\leq i,j\leq n}g^{ij}w_{,1ij}$: as we have
let
\begin{align}
w=&u-G,\quad G=\psi(x,u)d,\label{3w}
\end{align}
then we have
\begin{align}
w_1=&(1-G_u)u_1-G_{x_1},\notag\\
w_{,1i}=&(1-G_u)u_{,1i}-G_{uu}u_1u_i-G_{ux_i}u_1-G_{x_1u}u_i-G_{x_1x_i},\label{3wki}\\
w_{,1ij}=&(1-G_u)u_{,1ij}-G_{uu}(u_{,1i}u_j+u_{,1j}u_i+u_{,ij}u_1)\notag\\
&-G_{ux_i}u_{,1j}-G_{ux_j}u_{,1i}-G_{x_1u}u_{,ij}-G_{uuu}u_1u_iu_j\notag\\
&-G_{uux_j}u_1u_i-G_{ux_iu}u_1u_j-G_{x_1uu}u_iu_j\notag\\
&-G_{ux_ix_j}u_1-G_{x_1ux_j}u_i-G_{x_1x_iu}u_j-G_{x_1x_ix_j}.\label{3wkij}
\end{align}
So from the choice of the coordinate and the equations \eqref{equ}, \eqref{3wkij}, we have
\begin{align}
\sum_{1\leq i,j\leq n}g^{ij}w_{,1ij}
\ge&(1-G_u)\sum_{1\leq i,j\leq n}g^{ij}u_{,1ij}-2\sum_{1\leq i,j\leq n}g^{ij}(G_{uu}u_i+G_{ux_i})u_{,1j}\notag\\&-fG_{uu}vu_1-C_2u_1.\label{3wkaijwkij}
\end{align}

Differentiating \eqref{equ}, we have
\begin{align}\label{3aijuijka}
g^{ij}u_{,ij1}=&-\sum_{1\leq l\leq n}g^{ij}_{p_l}u_{,l1}u_{,ij}+vD_{1}f+fv_{1}.
\end{align}
From \eqref{2.4}, we have
\begin{align}\label{3aijpl}
g^{ij}_{p_l}=&-\frac{1}{v^2}(\delta_{il}u_j+\delta_{jl}u_i)+\frac{2}{v^4}u_iu_ju_l.
\end{align}
By the definition of $v$, we have
\begin{align}\label{3vvk}
v v_1=&\sum_{1\leq l\leq n} u_lu_{,l1}.
\end{align}
and
\begin{align}\label{3Dkf}
D_1f=&f_uu_1+f_{x_1}.
\end{align}
Inserting \eqref{3aijpl}, \eqref{3vvk} and \eqref{3Dkf} into \eqref{3aijuijka}, we have
\begin{align}
g^{ij}u_{,ij1}
=&\frac{2}{v^2}\sum_{1\leq l\leq n}g^{il}u_{,j}u_{,ij}u_{,l1}+f\sum_{1\leq l\leq n}\frac{u_lu_{,l1}}{v}+f_uvu_1+f_{x_1}v.\label{3aijuijkb1}
\end{align}
From the Ricci identity, we have
\begin{align}
u_{,ij1}
=&u_{,1ij}+\sum_{1\leq l\leq n}u_lR^l_{ij1}.\label{3aijuijkb2}
\end{align}
Inserting  \eqref{3aijuijkb2} into \eqref{3aijuijkb1},  we have
\begin{align}
g^{ij}u_{,1ij}
=&\frac{2}{v^2}\sum_{1\leq l\leq n}g^{il}u_{,j}u_{,ij}u_{,l1}-\sum_{1\leq l\leq n}u_lR^l_{ij1}+f\sum_{1\leq l\leq n}\frac{u_lu_{,l1}}{v}+f_uvu_1+f_{x_1}v.\label{3aijuijkb}
\end{align}
Inserting  \eqref{3aijuijkb} into \eqref{3wkaijwkij}, we have the formula for $I_{11}$
\begin{align}\label{3I11a}
I_{11}
\ge&\frac{2(1-G_u)w_1}{v^2}\sum_{1\leq i,j,l\leq n}g^{il}u_{,j}u_{,ij}u_{,l1}-2(1-G_u)w_1\sum_{1\leq i,j\leq n}g^{ij}(G_{uu}u_i+G_{ux_i})u_{,1j}
\notag\\&+f(1-G_u)w_1\sum_{1\leq l\leq n}\frac{u_lu_{,l1}}{v}-fG_{uu}(1-G_u)w_1vu_1-C_3u_1^2.
\end{align}
In the above inequality, we have used $f_u\ge 0$.

For the term $I_{12}$:
\begin{align}
I_{12}
\ge&\sum_{1\leq i,j\leq n}g^{ij}w_{1i}w_{1j}+\sum_{2\leq i\leq n}g^{ii}w^2_{ii}.\label{3I12a}
\end{align}

In the following we use the relation $\varphi_i(x_0)=0$, we get the formula
\begin{align}
w_{1i}=&-w_1\log w_1(h'u_i+g'd_i).\label{3w1ia}
\end{align}
From \eqref{3wki} and \eqref{3w1ia}, we have
\begin{align}
(1-G_u)u_{,1i}=&-w_1\log w_1(h'u_i+g'd_i)+[G_{uu}u_1u_i+G_{ux_i}u_1+G_{x_1u}u_i+G_{x_1x_i}].\label{3u1ia}
\end{align}
So putting \eqref{3wki},  \eqref{3w1ia}, \eqref{3u1ia} into \eqref{3I11a} and \eqref{3I12a}, we have
\begin{align}\label{3I11a3I12a}
I_{11}+I_{12}\geq&\frac{2(1-G_u)w_1}{v^2}\sum_{1\leq i,j,l\leq n}g^{il}u_{,j}u_{,ij}u_{,l1}+\sum_{1\leq i,j\leq n}g^{ij}w_{1i}w_{1j}+\sum_{2\leq i\leq n}g^{ii}w^2_{ii}\notag\\&
-2(1-G_u)w_1\sum_{1\leq i,j\leq n}g^{ij}(G_{uu}u_i+G_{ux_i})u_{,1j}
+f(1-G_u)w_1\sum_{1\leq l\leq n}\frac{u_lu_{,l1}}{v}\notag\\&-fG_{uu}(1-G_u)w_1vu_1-C_3u_1^2\notag\\
\geq &3h'^2w^2_1\log^2w_1\sum_{1\leq i,j\leq n}g^{ij}u_iu_j+3g'^2w^2_1\log^2w_1\sum_{1\leq i,j\leq n}g^{ij}d_id_j\notag\\&-h'fw^2_1\log^2w_1\frac{u_1^2}{v}-C_4w^2_1\log w_1\notag\\
=&3w^2_1\log^2w_1\sum_{1\leq i,j\leq n}g^{ij}u_iu_j+3\alpha_0^2w^2_1\log^2w_1\sum_{1\leq i,j\leq n}g^{ij}d_id_j\notag\\&-fw^2_1\log^2w_1\frac{u_1^2}{v}-C_4w^2_1\log w_1,
\end{align}
The above last equality has been used the choice of $h$ and $g$. It follows that from \eqref{3I1b},
\begin{align}\label{3I1c}
I_1\geq &3\log w_1\sum_{1\leq i,j\leq n}g^{ij}u_iu_j+3\alpha_0^2\log w_1\sum_{1\leq i,j\leq n}g^{ij}d_id_j-f\log w_1\frac{u_1^2}{v}-C_4.
\end{align}
Combining \eqref{3I1c}, \eqref{3I2c} and \eqref{3aijvarphiija}, it follows that
\begin{align}\label{3aijvarphiijb}
0\geq &\sum_{1\leq i,j\leq n}g^{ij}\varphi_{,ij}(x_0)\geq \frac{1}{4}\log w_1-C_5.
\end{align}
So there exists a positive constant
$C_{6}$ such that
 \begin{align}\label{3C13}
 |Du|(x_0)\leq C_{6}.
 \end{align}
So from  Case I, Case II,  and \eqref{3C13}, we have
 $$|Du|(x_0)\leq C_{7}, \quad \quad x_0\in\Omega_{\mu_0}\bigcup\partial\Omega.$$
Here the above $C_{2},\ldots, C_{7}$ are positive constants  depending  only
on $n, \Omega, \mu_0, M_0, \sup_{\overline{\Omega}}|Ric|, L_1, L_2$.

Since $\varphi(x)\leq\varphi(x_0),\quad \text{for} \quad x\in \Omega_{\mu_0}$, there exists $M_2$ such that
\begin{align}\label{3M2}
|Du|(x)\leq M_2, \quad in\quad\Omega_{\mu_0}\bigcup\partial\Omega,
\end{align}
where $M_2$ depends only on $n, \Omega, \sup_{\overline{\Omega}}|Ric|, \mu_0, M_0,  \sup_{\overline{\Omega}}|Ric|, L_1, L_2$.

So at last we get the following estimate
$$\sup_{\overline\Omega_{\mu_0}}|Du|\leq \max\{M_1, M_2\},$$
where the positive constant  $ M_1$ depends only on $n, \mu_0, M_0,\sup_{\overline{\Omega}}|Ric|, L_1$; and $ M_2$ depends only on ~$n, \Omega, \mu_0, M_0, \sup_{\overline{\Omega}}|Ric|, L_1, L_2$.

 So  we  complete the proof of Theorem~\ref{Thm1.1}.\qed

As a consequence, from the standard $C^0$  estimates for the prescribed mean curvature equation in Concus-Finn \cite{CF74} ( see also Spruck \cite{Sp75}, Ma-Xu\cite{MX14}), we can also get the  existence theorem Theorem~\ref{Thm1.2}.

 {\bf Acknowledgement.}
The authors would like to
thank Professor  Xinan Ma for  suggesting this question and  helpful discussions on this paper.

\end{document}